\documentclass[11pt]{amsart}

\usepackage{amsmath}
\usepackage{amssymb}
\usepackage{mathrsfs}
\usepackage[tight]{subfigure}
\usepackage{epsfig,epstopdf,psfrag}
\usepackage{array,mathtools}
\usepackage{booktabs}
\usepackage{natbib}

\newtheorem{theorem}{Theorem}

\newtheorem{definition}[theorem]{Definition}
\newenvironment{remark}{\noindent \textbf{Remark}.}{\hfill $\square$}
\newcommand{\divg}{\mathop{\rm div}\nolimits}
\newcommand{\grad}{\mathop{\rm grad}\nolimits}

\providecommand{\abs}[1]{\lvert#1\rvert}
\providecommand{\norm}[1]{\lvert#1\rvert}

\numberwithin{equation}{section}
\numberwithin{theorem}{section}
\numberwithin{figure}{subsection}
\numberwithin{table}{subsection}

\title{Pointwise Inequalities for Elliptic Boundary Value Problems}
\author{Guo Luo \and Vladimir G. Maz'ya}
\begin{document}
\maketitle

\section*{Abstract}
We introduce a new approach to obtaining pointwise estimates for solutions of elliptic boundary value problems when the operator being considered satisfies a certain type of weighted integral inequalities. The method is illustrated on several examples, including a scalar second-order elliptic equation, the 3D Lam\'{e} system, and a scalar higher-order elliptic equation. The techniques can be extended to other elliptic boundary value problems provided that the corresponding weighted integral inequalities are satisfied.
%
%


\section{Introduction}\label{sec_intro}
An important open problem in the mathematical theory of linear elasticity is whether solutions of the elasticity system, when supplemented with homogeneous Dirichlet boundary conditions and sufficiently smooth right-hand side data, are uniformly bounded in \emph{arbitrary} domains. A similar question stands in the theory of hydrostatics, where the uniform boundedness of solutions of the Stokes system in general domains remains unknown. For bounded domains $\Omega$ with smooth boundaries $\partial \Omega$, inequalities of the form
\begin{equation}
  \norm{u}_{L^{\infty}(\Omega)} \leq C_{\Omega} \norm{D^{k} u}_{L^{q}(\Omega)}^{a} \norm{Lu}_{W^{l,p}(\Omega)}^{b}
  \label{eqn_pt_bd}
\end{equation}
can often be obtained, by combining appropriate \emph{a priori} estimates with Sobolev inequalities. The problem of such inequalities is that the constant $C_{\Omega}$ generally depends on the smoothness of the domain $\Omega$, and can blow up if the boundary of $\Omega$ contains geometric singularities.

Efforts have been devoted to the study of inequalities of the type \eqref{eqn_pt_bd} with constants \emph{independent} of the domain $\Omega$. In \citet{xie1991}, a sharp pointwise bound
\begin{equation}
  \norm{u}_{L^{\infty}(\Omega)}^{2} \leq \frac{1}{2\pi}\, \norm{Du}_{L^{2}(\Omega)} \norm{\Delta u}_{L^{2}(\Omega)}
  \label{eqn_pt_bd_lap}
\end{equation}
was obtained for functions $u$ with zero-trace and with $L^{2}$-integrable gradient $Du$ and Laplacian $\Delta u$ on arbitrary three-dimensional domains $\Omega$. The constant $(2\pi)^{-1}$ was shown to be the best possible. In \citet{xie1992}, a similar inequality with a slightly different best constant ($(3\pi)^{-1}$ instead of $(2\pi)^{-1}$) was conjectured for the Stokes system on arbitrary domains $\Omega \subseteq \mathbb{R}^{3}$, and was proved in the special case $\Omega = \mathbb{R}^{3}$. Further development and results along these lines can be found in \citet{heywood2001} and the references there in. Regarding the 3D Lam\'{e} system, estimates of the form \eqref{eqn_pt_bd_lap} seem to be less well studied, and we are not aware of any results or conjectures similar to \eqref{eqn_pt_bd_lap}.

In this paper, we introduce a new approach to obtaining pointwise inequalities of the form \eqref{eqn_pt_bd} with constants independent of the domain $\Omega$. The method works for elliptic operators $L$ satisfying a weighted integral inequality
\begin{equation}
  \int_{\Omega} Lu \cdot \Phi u\,dx \geq 0,
  \label{eqn_wpd}
\end{equation}
where the weight $\Phi$ is either the fundamental solution or Green's function of $L$. Weighted inequalities of the form \eqref{eqn_wpd} were first established for second-order scalar operators (see Section \ref{ssec_pt_bd_1} below for a prototypical derivation), and were later generalized to certain higher-order scalar operators \citep{mazya2002} and second-order systems \citep{lm2007,lm2010}. They have important applications in the regularity theory of boundary points, and have been studied extensively in \citet{mazya1977,mazya1979,mazya1999,mazya2002,mazya1991,eilertsen2000}. By utilizing a slightly modified version of \eqref{eqn_wpd} (see \eqref{eqn_wpd_sx} below), we shall show that the pointwise estimate \eqref{eqn_pt_bd} follows almost immediately from the weighted positivity of $L$, and the constant thus obtained is independent of the domain $\Omega$. The method is illustrated on several concrete problems, including a scalar second-order elliptic equation, the 3D Lam\'{e} system, and a scalar higher-order elliptic equation. The techniques can be extended to other elliptic boundary value problems provided that the corresponding weighted integral inequalities are satisfied.

In what follows, $W_{0}^{k,q}(\Omega)$ denotes the usual Sobolev space of zero-trace functions, equipped with the norm
\begin{displaymath}
  \norm{u}_{W_{0}^{k,q}(\Omega)} := \biggl( \sum_{\abs{\beta} \leq k} \norm{D^{\beta} u}_{L^{q}(\Omega)}^{q} \biggr)^{1/q}.
\end{displaymath}

\begin{theorem}
Let $L$ be a second-order elliptic operator,
\begin{displaymath}
  Lu = -D_{i} (a_{ij}(x) D_{j} u),\qquad D_{i} = \frac{\partial}{\partial x_{i}},
\end{displaymath}
defined in a bounded domain $\Omega \subset \mathbb{R}^{n}\ (n \geq 3)$ with real, measurable coefficients $a_{ij}(x)$. Suppose $L$ satisfies the strong ellipticity condition
\begin{displaymath}
  \lambda \abs{\xi}^{2} \leq a_{ij}(x) \xi_{i} \xi_{j} \leq \Lambda \abs{\xi}^{2},\qquad \Lambda \geq \lambda > 0,
\end{displaymath}
for almost all $x \in \Omega$ and all $\xi = (\xi_{1},\dotsc,\xi_{n}) \in \mathbb{R}^{n}$. Let $s < n/(n-2),\ p = s/(s-1),\ q = (n-2)s$ and let $u \in W_{0}^{1,q}(\Omega)$ be such that $Lu \in L^{p}(\Omega)$. Then
\begin{displaymath}
  \norm{u}_{L^{\infty}(\Omega)} \leq C \norm{Lu}_{L^{p}(\Omega)}^{a} \norm{Du}_{L^{q}(\Omega)}^{1-a},\qquad a = \frac{1}{n-1},
\end{displaymath}
where $C$ is an absolute constant depending only on $\lambda,\ \Lambda,\ n$, and $s$.
\label{thm_pt_bd_1}
\end{theorem}

\begin{theorem}
Let $L$ be the 3D Lam\'{e} system,
\begin{equation}
  Lu = -\Delta u - \alpha \grad \divg u = -D_{kk} u_{i} - \alpha D_{ki} u_{k},\qquad i = 1,2,3,
  \label{eqn_lame}
\end{equation}
where $\alpha = 1/(1-2\nu) > -1$ and $\nu$ is Poisson's ratio. Let $\alpha \in (\alpha_{-},\alpha_{+}) \approx (-0.194,1.524),\ q < 3,\ p = q/(q-1)$ and let $u \in W_{0}^{1,q}(\Omega)$ be such that $Lu \in L^{p}(\Omega)$, where $\Omega$ is an arbitrary bounded domain in $\mathbb{R}^{3}$. Then
\begin{displaymath}
  \norm{u}_{L^{\infty}(\Omega)} \leq C \norm{Lu}_{L^{p}(\Omega)}^{1/2} \norm{Du}_{L^{q}(\Omega)}^{1/2},
\end{displaymath}
where $C$ is an absolute constant depending only on $\alpha$ and $q$.
\label{thm_pt_bd_lame}
\end{theorem}

\begin{theorem}
Let $L$ be an elliptic operator of order $2m$,
\begin{displaymath}
  L = (-1)^{m} \sum_{\abs{\alpha} = \abs{\beta} = m} a_{\alpha \beta} D^{\alpha+\beta},\qquad D^{\alpha} = \frac{\partial^{\abs{\alpha}}}{\partial x_{1}^{\alpha_{1}} \partial x_{2}^{\alpha_{2}} \dotsb \partial x_{n}^{\alpha_{n}}},
\end{displaymath}
defined in $\mathbb{R}^{n}\ (n > 2m)$ with real, constant coefficients $a_{\alpha\beta}$. Suppose $L$ satisfies the strong ellipticity condition
\begin{displaymath}
  \lambda \abs{\xi}^{2m} \leq a_{\alpha\beta} \xi^{\alpha} \xi^{\beta} \leq \Lambda \abs{\xi}^{2m},\qquad \Lambda \geq \lambda > 0,
\end{displaymath}
for all $\xi = (\xi_{1},\dotsc,\xi_{n}) \in \mathbb{R}^{n}$. Let $F$ be the fundamental solution of $L$. Let $q < n/(n-2m),\ q' = q/(q-1),\ k = n-2m$ and let $u \in W_{0}^{k,q}(\Omega)$ be such that $Lu \in L^{q'}(\Omega)$, where $\Omega$ is an arbitrary bounded domain in $\mathbb{R}^{n}$. If $F$ is homogeneous of order $2m-n$ and $L$ is weighted positive with the weight $F$, then
\begin{displaymath}
  \norm{u}_{L^{\infty}(\Omega)}^{2} \leq C \norm{D^{k} u}_{L^{q}(\Omega)} \norm{Lu}_{L^{q'}(\Omega)},
\end{displaymath}
where $C$ is an absolute constant depending only on $\lambda,\ \Lambda,\ m,\ n$ and $q$.
\label{thm_pt_bd_m}
\end{theorem}

The proofs of Theorems \ref{thm_pt_bd_1}--\ref{thm_pt_bd_m} are given in Section \ref{sec_pt_bd}.

\section{The Notion of Weighted Positivity}\label{sec_wpd}
Let $\Omega$ be a domain in $\mathbb{R}^{n}$ and let $L = (L_{i})_{i=1}^{N}$ be a scalar elliptic operator ($N = 1$) or an elliptic system ($N > 1$) defined on $\Omega$. Without making further structural assumptions on $L$, we recall first the abstract notion of weighted positivity.

\begin{definition}[Weighted positivity]
Assume that $0 \in \Omega$ and that $\Psi(x) = (\Psi_{ij}(x))_{i,j=1}^{N}$ is a given (matrix) function that is sufficiently regular except possibly at $x = 0$. The operator $L = (L_{i})_{i=1}^{N}$ is said to be weighted positive ($N = 1$) or weighted positive definite ($N > 1$) with the weight $\Psi$ if
\begin{equation}
  \int_{\Omega} Lu \cdot \Psi u\,dx = \int_{\Omega} (Lu)_{i} \Psi_{ij} u_{j}\,dx \geq 0,
  \label{eqn_wpd_w}
\end{equation}
for all real-valued, smooth vector functions $u = (u_{i})_{i=1}^{N} \in C_{0}^{\infty}(\Omega \setminus \{0\})$.
\end{definition}

The concept of weighted positivity is often more useful when the integral $\int_{\Omega} Lu \cdot \Psi u\,dx$ in \eqref{eqn_wpd_w} has a \emph{positive} lower bound. This motivates the following definition.

\begin{definition}[Strong weighted positivity]
Assume that $0 \in \Omega,\ L = (L_{i})_{i=1}^{N}$ is of order $2m,\ m \geq 1$, and that $\Psi(x) = (\Psi_{ij}(x))_{i,j=1}^{N}$ is a given (matrix) function that is sufficiently regular except possibly at $x = 0$. The operator $L$ is said to be strongly weighted positive ($N = 1$) or strongly weighted positive definite ($N > 1$) with the weight $\Psi$ if, for some $c > 0$,
\begin{equation}
  \int_{\Omega} Lu \cdot \Psi u\,dx \geq c \sum_{k=1}^{m} \int_{\Omega} \abs{D^{k} u}^{2} \abs{x}^{2k-2m} \abs{\Psi}\,dx,
  \label{eqn_wpd_s}
\end{equation}
for all real-valued, smooth vector functions $u = (u_{i})_{i=1}^{N} \in C_{0}^{\infty}(\Omega \setminus \{0\})$. Here $D^{k}$ denotes the gradient operator of order $k$, i.e. $D^{k} = \{D^{\alpha}\}$ with $\abs{\alpha} = k$, and $\abs{\Psi}$ stands for the Frobenius norm of $\Psi$, i.e. $\abs{\Psi}^{2} = \sum_{i,j=1}^{N} \abs{\Psi_{ij}}^{2}$.
\end{definition}

Among all possible candidates of the weight function $\Psi$, the special choice $\Psi = \Phi$ where $\Phi$ is the fundamental solution or Green's function of $L$ is of the most interest. In particular, if $L$ has constant coefficients and is strongly weighted positive in the sense of \eqref{eqn_wpd_s}, with the weight $\Psi = \Phi$, then it can be shown that, for all $x \in \Omega$ and for the same constant $c$ as given in \eqref{eqn_wpd_s},
\begin{equation}
  \int_{\Omega} Lu \cdot \Phi(x-y) u\,dy \geq \frac{1}{2}\, \abs{u(x)}^{2} + c \sum_{k=1}^{m} \int_{\Omega} \frac{\abs{D^{k} u}^{2}}{\abs{x-y}^{2m-2k}}\, \abs{\Phi(x-y)}\,dy,
  \label{eqn_wpd_sx}
\end{equation}
for all real-valued, smooth vector functions $u = (u_{i})_{i=1}^{N} \in C_{0}^{\infty}(\Omega)$ \citep[see, for example,][]{mazya2002}. Estimate \eqref{eqn_wpd_sx} is significant because it provides a pointwise bound for the test function $u$, and it serves as the basis of the estimates to be derived below in Section \ref{sec_pt_bd}.

\section{Proofs of the Main Theorems}\label{sec_pt_bd}
\subsection{Proof of Theorem \ref{thm_pt_bd_1}}\label{ssec_pt_bd_1}
We shall show that the multiplicative inequality
\begin{equation}
  \norm{u}_{L^{\infty}(\Omega)}^{n-1} \leq C \norm{Lu}_{L^{p}(\Omega)} \norm{Du}_{L^{q}(\Omega)}^{n-2},\qquad p = \frac{s}{s-1},\ q = (n-2)s,
  \label{eqn_pt_bd_1}
\end{equation}
holds for all $s < n/(n-2)$ with a constant $C = C(\lambda,\Lambda,n,s)$ independent of the domain $\Omega$, if $u \in W_{0}^{1,q}(\Omega)$ and $Lu \in L^{p}(\Omega)$. Inequality \eqref{eqn_pt_bd_1} will be obtained by a modification of weighted positivity of the operator $L$ with zero Dirichlet boundary conditions. More precisely, let $G(x,y)$ denote Green's function of $L$. Its existence and uniqueness are classical facts as well as the estimates
\begin{equation}
  c_{1} \abs{x-y}^{2-n} \leq G(x,y) \leq c_{2} \abs{x-y}^{2-n},
  \label{eqn_G_1}
\end{equation}
where $c_{1}$ and $c_{2}$ are positive constants depending on $\lambda$ and $\Lambda$ \citep[see][]{royden1962,lsw1963}. By definition of a weak solution of the Dirichlet problem and by a standard approximation argument, we obtain, for almost all $x \in \Omega$,
\begin{multline*}
  \int_{\Omega} Lu \cdot G(x,y) u \abs{u}^{n-3}\,dy \\
  = (n-2) \int_{\Omega} a_{ij} D_{i} u D_{j} u \cdot G(x,y) \abs{u}^{n-3}\,dy + \int_{\Omega} a_{ij} D_{y_{i}} G(x,y) D_{j} u \cdot u \abs{u}^{n-3}\,dy \\
  \geq \frac{1}{n-1} \int_{\Omega} a_{ij} D_{y_{i}} G(x,y) D_{j} \abs{u}^{n-1}\,dy = \frac{1}{n-1}\, \abs{u(x)}^{n-1}.
\end{multline*}
Hence for $s < n/(n-2)$,
\begin{displaymath}
  \abs{u(x)}^{n-1} \leq c_{2} (n-1) \norm{Lu}_{L^{p}(\Omega)} \biggl( \int_{\Omega} \frac{\abs{u}^{q}}{\abs{x-y}^{q}}\,dy \biggr)^{1/s},\qquad p = \frac{s}{s-1},\ q = (n-2)s,
\end{displaymath}
where $c_{2}$ is the constant in \eqref{eqn_G_1}. By Hardy's inequality, the last integral does not exceed
\begin{displaymath}
  \biggl( \frac{q}{n-q} \biggr)^{n-2} \biggl( \int_{\Omega} \abs{Du}^{q}\,dy \biggr)^{1/s},
\end{displaymath}
thus \eqref{eqn_pt_bd_1} follows with
\begin{displaymath}
  C = c_{2} (n-1) \biggl( \frac{q}{n-q} \biggr)^{n-2}.
\end{displaymath}
This completes the proof of Theorem \ref{thm_pt_bd_1}.

\begin{remark}
Note that for $L = -\Delta$, we have $c_{2} = [(n-2) \omega_{n}]^{-1}$ where $\omega_{n}$ is the measure of the unit $(n-1)$-sphere $S^{n-1}$. In particular, for $n = 3,\ s = 2$, and $L = -\Delta$, we have $p = 2,\ q = 2$, and $c_{2} = (4\pi)^{-1}$. Thus Theorem \ref{thm_pt_bd_1} implies that
\begin{displaymath}
  \pi \norm{u}_{L^{\infty}(\Omega)}^{2} \leq \norm{\Delta u}_{L^{2}(\Omega)} \norm{Du}_{L^{2}(\Omega)}.
\end{displaymath}
This is similar to the pointwise bound \eqref{eqn_pt_bd_lap} obtained by \citet{xie1991}, but with a slightly worse constant ($\pi$ instead of $2\pi$).
\end{remark}

\begin{remark}
Note that the inequality \eqref{eqn_pt_bd_1} fails for $s = n/(n-2)$. Indeed, when $s = n/(n-2)$, it is easily checked that $p = n/2$ and $q = n$. Let $\zeta \in C_{0}^{\infty}[0,\infty)$ be a smooth cutoff function on $\mathbb{R}$ such that $0 \leq \zeta \leq 1,\ \zeta(x) = 1$ for $x \leq 1/2$, and $\zeta = 0$ for $x \geq 1$. Set
\begin{displaymath}
  u(x) = \zeta(\abs{x}) \log \abs{\log\abs{x}}.
\end{displaymath}
It is easily verified that $u$ is unbounded in any neighborhood of $x = 0$ while
\begin{displaymath}
  \abs{Du(x)} \leq C \abs{x}^{-1} \abs{\log\abs{x}}^{-1},\qquad \abs{\Delta u(x)} \leq C \abs{x}^{-2} \abs{\log\abs{x}}^{-1},
\end{displaymath}
for small $\abs{x}$, indicating that $Du \in L^{n}(\mathbb{R}^{n})$ and $\Delta u \in L^{n/2}(\mathbb{R}^{n})$. This violates \eqref{eqn_pt_bd_1}.
\end{remark}

\subsection{Proof of Theorem \ref{thm_pt_bd_lame}}\label{ssec_pt_bd_lame}
Let $\Phi$ be the fundamental matrix of the 3D Lam\'{e} operator $L$,
\begin{subequations}\label{eqn_lame_fs}
\begin{align}
  \Phi_{ij}(x) & = c_{\alpha} r^{-1} \biggl( \delta_{ij} + \frac{\alpha}{\alpha+2}\, \omega_{i} \omega_{j} \biggr),\qquad i,j = 1,2,3, \label{eqn_lame_fs_phi} \\
  c_{\alpha} & = \frac{\alpha+2}{8\pi(\alpha+1)} > 0, \label{eqn_lame_fs_c}
\end{align}
\end{subequations}
where $\delta_{ij}$ is the Kronecker delta, $r = \abs{x}$, and $\omega_{i} = x_{i}/\abs{x}$. The weighted positive definiteness of $L$ with the weight $\Phi$ has been established in \citet{lm2010} for certain ranges of the parameter $\alpha$.

\begin{theorem}[\citealp{lm2010}]
The 3D Lam\'{e} system $L$ is weighted positive definite with the weight $\Phi$ when $\alpha_{-} < \alpha < \alpha_{+}$, where $\alpha_{-} \approx -0.194$ and $\alpha_{+} \approx 1.524$. It is not weighted positive definite with the weight $\Phi$ when $\alpha < \alpha_{-}^{(c)} \approx -0.902$ or $\alpha > \alpha_{+}^{(c)} \approx 39.450$.
\label{thm_lame_wpd}
\end{theorem}

By examining the proof of Theorem \ref{thm_lame_wpd}, it can be observed that $L$ is in fact \emph{strongly} weighted positive definite with the weight $\Phi$, i.e.
\begin{equation}
  \int_{\mathbb{R}^{3}} Lu \cdot \Phi u\,dx \geq \frac{1}{2}\, \abs{u(0)}^{2} + c \int_{\mathbb{R}^{3}} \abs{Du}^{2} \abs{x}^{-1}\,dx,\qquad \text{for some $c > 0$},
  \label{eqn_wpd_lame}
\end{equation}
for all $u = (u_{i})_{i=1}^{3} \in C_{0}^{\infty}(\mathbb{R}^{3})$ when $\alpha_{-} < \alpha < \alpha_{+}$.

Using Theorem \ref{thm_lame_wpd}, we shall show that the multiplicative inequality
\begin{equation}
  \norm{u}_{L^{\infty}(\Omega)}^{2} \leq C \norm{Lu}_{L^{p}(\Omega)} \norm{Du}_{L^{q}(\Omega)},\qquad p = \frac{q}{q-1},
  \label{eqn_pt_bd_lame}
\end{equation}
holds on arbitrary bounded domains $\Omega \subset \mathbb{R}^{3}$ with a constant $C = C(\alpha,q)$ independent of the domain, provided that $\alpha_{-} < \alpha < \alpha_{+},\ q < 3,\ u \in W_{0}^{1,q}(\Omega)$ and $Lu \in L^{p}(\Omega)$. Inequality \eqref{eqn_pt_bd_lame} will be obtained as an immediate corollary of the weighted positive definiteness of $L$. More precisely, the weighted inequality \eqref{eqn_wpd_lame} implies that
\begin{displaymath}
  \int_{\mathbb{R}^{3}} Lu \cdot \Phi(x-y) u\,dy \geq \frac{1}{2}\, \abs{u(x)}^{2},
\end{displaymath}
for all real-valued $u = (u_{i})_{i=1}^{3} \in C_{0}^{\infty}(\mathbb{R}^{3})$. By definition of a weak solution and by a standard approximation argument, the same inequality can be proved for $u \in W_{0}^{1,q}(\Omega)$ for which $Lu \in L^{p}(\Omega)$. Hence for all such $u$ and all $q < 3$, we have
\begin{displaymath}
  \abs{u(x)}^{2} \leq 2c_{3} \norm{Lu}_{L^{p}(\Omega)} \biggl( \int_{\Omega} \frac{\abs{u}^{q}}{\abs{x-y}^{q}}\,dy \biggr)^{1/q},\qquad c_{3} = c_{\alpha} \biggl( 1 + \frac{\abs{\alpha}}{\alpha+2} \biggr),\ p = \frac{q}{q-1},
\end{displaymath}
where $c_{\alpha}$ is the constant in \eqref{eqn_lame_fs_c}. By Hardy's inequality, the last integral does not exceed
\begin{displaymath}
  \frac{q}{3-q} \biggl( \int_{\Omega} \abs{Du}^{q}\,dy \biggr)^{1/q},
\end{displaymath}
thus \eqref{eqn_pt_bd_lame} follows with
\begin{displaymath}
  C = 2 c_{\alpha} \biggl( 1 + \frac{\abs{\alpha}}{\alpha+2} \biggr) \biggl( \frac{q}{3-q} \biggr).
\end{displaymath}
This completes the proof of Theorem \ref{thm_pt_bd_lame}.

\begin{remark}
The hydrostatic limit $\alpha \to \infty$ of the 3D Lam\'{e} system lies unfortunately outside the regime of weighted positive definiteness, hence the uniform estimates of solutions of the 3D Stokes system, if they exist, cannot be deduced from the weighted inequality \eqref{eqn_wpd_lame}.
\end{remark}

\subsection{Proof of Theorem \ref{thm_pt_bd_m}}\label{ssec_pt_bd_m}
Let $F(x)$ denote the fundamental solution of $L$. It is well known that $F$ exists for all $n > 2m$ and is homogeneous of order $2m-n$,
\begin{equation}
  F(x) = \abs{x}^{2m-n} F \biggl( \frac{x}{\abs{x}} \biggr),\qquad x \in \mathbb{R}^{n} \setminus \{0\},
  \label{eqn_F_m}
\end{equation}
when $n$ is odd \citep{john1982}. When $n$ is even, \eqref{eqn_F_m} may not be valid since terms of the order $\abs{x}^{2m-n} \log\abs{x}$ may occur in $F$. Under the assumptions that $L$ is weighted positive with the weight $F$, i.e.
\begin{equation}
  \int_{\mathbb{R}^{n}} Lu \cdot F u\,dx \geq 0,
  \label{eqn_wpd_m}
\end{equation}
for all real-valued $u \in C_{0}^{\infty}(\mathbb{R}^{n} \setminus \{0\})$, and that $F$ satisfies \eqref{eqn_F_m}, we shall show that the multiplicative inequality
\begin{equation}
  \norm{u}_{L^{\infty}(\Omega)}^{2} \leq C \norm{D^{k} u}_{L^{q}(\Omega)} \norm{Lu}_{L^{q'}(\Omega)},\qquad k = n-2m,\ q' = \frac{q}{q-1},
  \label{eqn_pt_bd_m}
\end{equation}
holds on arbitrary bounded domains $\Omega \subset \mathbb{R}^{n}$ with a constant $C = C(\lambda,\Lambda,n,m,q)$ independent of the domain, provided that $q < n/(n-2m),\ u \in W_{0}^{k,q}(\Omega)$ and $Lu \in L^{q'}(\Omega)$. Like in the previous two theorems, inequality \eqref{eqn_pt_bd_m} will be obtained as an immediate corollary of the weighted positivity of $L$. More precisely, the weighted inequality \eqref{eqn_wpd_m} implies that
\begin{displaymath}
  \int_{\mathbb{R}^{n}} Lu \cdot F(x-y) u\,dy \geq \frac{1}{2}\, \abs{u(x)}^{2},
\end{displaymath}
for all real-valued $u \in C_{0}^{\infty}(\mathbb{R}^{n})$ \citep[see][Proposition 3]{mazya2002}. By definition of a weak solution and by a standard approximation argument, the same inequality can be proved for $u \in W_{0}^{k,q}(\Omega)$ for which $Lu \in L^{q'}(\Omega)$. Hence for all such $u$ and all $q < n/(n-2m)$, we have
\begin{displaymath}
  \abs{u(x)}^{2} \leq 2c_{4} \norm{Lu}_{L^{q'}(\Omega)} \biggl( \int_{\Omega} \frac{\abs{u}^{q}}{\abs{x-y}^{kq}}\,dy \biggr)^{1/q},\qquad k = n-2m,\ q' = \frac{q}{q-1},
\end{displaymath}
where $c_{4} = \max_{\omega \in S^{n-1}} \abs{F(\omega)}$ is a positive constant depending on $\lambda$ and $\Lambda$. By repeated applications of Hardy's inequality, the last integral does not exceed
\begin{displaymath}
  \biggl( \frac{1}{r-k} \biggr) \biggl( \frac{1}{r-k+1} \biggr) \dotsb \biggl( \frac{1}{r-1} \biggr) \biggl( \int_{\Omega} \abs{D^{k} u}^{q}\,dy \biggr)^{1/q},\qquad r = n/q,
\end{displaymath}
thus \eqref{eqn_pt_bd_m} follows with
\begin{displaymath}
  C = 2c_{4} \biggl( \frac{1}{r-k} \biggr) \biggl( \frac{1}{r-k+1} \biggr) \dotsb \biggl( \frac{1}{r-1} \biggr).
\end{displaymath}
This completes the proof of Theorem \ref{thm_pt_bd_m}.

\begin{remark}
Note that for $L = (-\Delta)^{m}$ with $2m < n$, \eqref{eqn_wpd_m} is satisfied if and only if $n = 5,\ 6,\ 7$ for $m = 2$ \citep{mazya1979} and $n = 2m+1,\ 2m+2$ for $m > 2$ \citep{mazya1997}. In particular, for $m = 2$ and $q = 2$, we have $q' = 2$ and $c_{2} = [2(n-2)(n-4)\omega_{n}]^{-1}$ where
\begin{displaymath}
  \omega_{n} = \frac{n \pi^{n/2}}{\Gamma(\frac{1}{2} n+1)}
\end{displaymath}
is the measure of the unit $(n-1)$-sphere $S^{n-1}$. Thus
\begin{displaymath}
  \norm{u}_{L^{\infty}(\Omega)}^{2} \leq \frac{\Gamma(4-\frac{1}{2} n)}{2\pi^{n/2} (n-2)(n-4)}\, \norm{D^{n-4} u}_{L^{2}(\Omega)} \norm{\Delta^{2} u}_{L^{2}(\Omega)},\qquad n = 5,\ 6,\ 7.
\end{displaymath}
To the best of our knowledge, inequalities of this type have not been known before.
\end{remark}

\end{document}